\input amstex

\documentstyle{amsppt}

\refstyle{A}

\nologo

%\magnification=\magstep1 

\hoffset .25 true in
\voffset .2 true in

\hsize=6.1 true in
\vsize=8.2 true in

\define\Adot{\bold A^\bullet}

\define\Pdot{\bold P^\bullet}

\define\ic{IC^\bullet_{{}_X}}

\define\cx{k^\bullet_{{}_X}[n]}

\define\cu{k^\bullet_{{}_{\Cal U}}[n]}

\define\dm{\operatorname{dim}}

\define\ob{\overset\circ\to{B_\epsilon}}

\topmatter

\title Intersection Cohomology, Monodromy, and the Milnor Fiber \endtitle

\author David B. Massey \endauthor

\address{David B. Massey, Dept. of Mathematics, Northeastern University, Boston,
MA, 02115, USA} \endaddress

\email{dmassey\@neu.edu}\endemail

\abstract   We say that a complex analytic space, $X$, is an intersection cohomology manifold if and only if the shifted constant sheaf on $X$ is isomorphic to intersection cohomology; this is quickly seen to be equivalent to $X$ being a homology manifold. Given an analytic function $f$ on an intersection cohomology manifold, we describe a simple relation between $V(f)$ being an intersection cohomology manifold and the vanishing cycle Milnor monodromy of $f$. We then describe how the Sebastiani-Thom isomorphism allows us to easily produce intersection cohomology manifolds with arbitrary singular sets. Finally, as an easy application, we obtain restrictions on the cohomology of the Milnor fiber of a hypersurface with a special type of one-dimensional critical locus. \endabstract

\endtopmatter

\document

\noindent\S0. {\bf Introduction}  

\vskip .1in

Suppose that $X$ is a connected open neighborhood of the origin in $\Bbb C^n$, where $n\geqslant 2$, and suppose that $f:X\rightarrow\Bbb C$ is (complex) analytic. Let $\bold p\in V(f)$ and let $F_{f, \bold p}$ denote the Milnor fiber of $f$ at $\bold p$.

\vskip .1in

In Theorem 8.5 of Milnor's classic ``Singular Points of Complex Hypersurfaces'' [{\bf Mi}], Milnor proved that, if $n\neq 3$ and $f$ has an isolated critical point at $\bold p$, then the real link of $V(f)$ at $\bold p$ is a topological sphere if and only if the identity minus the (Milnor) monodromy, acting on $H^{n-1}(F_{f, \bold p};\ \Bbb Z)$,  is invertible over $\Bbb Z$.

\vskip .1in

Milnor's proof can easily be used in the case where $f$  may have non-isolated critical points; in this setting, Milnor's proof shows that $V(f)$ is an integral cohomology manifold (or integral homology manifold) if and only if, for all $\bold p\in V(f)$, the identity minus the monodromy, $\operatorname{id}-\widetilde T_{f, \bold p}$, acting on the reduced cohomology of the Milnor fiber of $f$ at $\bold p$, $\widetilde H^*(F_{f, \bold p}; \ k)$, induces isomorphisms (in all degrees). This was demonstrated in 1977 by Randell [{\bf R}]. Our version of this result appears in Theorem 2.2.

\vskip .1in

Let us work now over a base ring $k$ which is actually a field (as opposed to working over $\Bbb Z$).
The categorical implication of the previous paragraph is that $V(f)$ is a $k$-cohomology manifold if and only if the morphism $\operatorname{id}-\widetilde T_{f}$, acting on the vanishing cycles $\phi_f[-1]k_{{}_X}^\bullet[n]$, is an isomorphism in the derived category $D^b_c(V(f))$ and, hence, an isomorphism in the full subcategory $Perv(V(f))$ of perverse sheaves on $V(f)$. This leads one to ask: is saying that an analytic space $Y$ is a cohomology or homology manifold equivalent to some nice statement  in the category $Perv(Y)$?

As Borho and MacPherson showed in [{\bf B-M}] and as we recall in Theorem 1.1, the answer to the above question, is ``yes''. In fact, the axiomatic characterization of intersection cohomology in terms of the support and co-support conditions implies immediately that an $m$-dimensional complex analytic space $Y$ is a $k$-cohomology manifold if and only if the shifted constant sheaf $k_{{}_Y}^\bullet[m]$ is isomorphic to the (middle perversity) intersection cohomology complex $IC^\bullet_{{}_Y}$. Here, by intersection cohomology, we mean that we are using the indexing which places possibly non-zero cohomology only in non-negative degrees, i.e., if $Y$ were a purely $m$-dimensional {\bf manifold}, then $IC^\bullet_{{}_Y}$ would be isomorphic to the shifted constant sheaf $k_{{}_Y}^\bullet[m]$.

What is the importance of $k_{{}_Y}^\bullet[m]$ being isomorphic to $IC^\bullet_{{}_Y}$? It is of great categorical importance, for if $Y$ is irreducible, then $IC^\bullet_{{}_Y}$ is a simple object in the Abelian, locally Artinian category $Perv(Y)$. Thus, if $Y$ is irreducible, $Y$ is a $k$-cohomology manifold if and only if $k_{{}_Y}^\bullet[m]$ is a simple object in $Perv(Y)$. As it is this categorical property of cohomology manifolds that we make use of throughout our arguments, we prefer to refer to such spaces $Y$ not as $k$-cohomology manifolds, but, rather,  as {\it $k$-intersection cohomology manifolds} (or $k$-$IC^\bullet$ manifolds).

\vskip .1in

In Theorem 2.2, we demonstrate the usefulness of the categorical approach. We consider $f:X\rightarrow\Bbb C$, but no longer require $X$ to be smooth -- instead, we allow $X$ itself to be a $k$-$IC^\bullet$ manifold. In this generalized setting, we relate $V(f)$ being a $k$-$IC^\bullet$ manifold to $\operatorname{id}-\widetilde T_{f}$ being an isomorphism in $Perv(V(f))$. The proof that we give is very short (though significant knowledge of the derived category is required); moreover, the hypotheses of the theorem are as weak as possible and are natural from the categorical point-of-view. Much of Theorem 2.2 could be proved along the lines of Milnor's argument, but the argument would be longer, special cases would be problematic, and the optimal hypotheses would be difficult to arrive at.

\vskip .2in

Let us consider a simple, but not trivial, example: $k=\Bbb C$, $X=\Bbb C^2$, and $f:X\rightarrow\Bbb C$ is given by $f(xy)=xy$. Let $m$ denote the inclusion of the origin, $\bold 0$, into $V(f)$. Let $j^x$ (respectively, $j^y$) denote the inclusion of $V(y)$ (respectively, $V(x)$) into $V(f)$.

It is well-known that the reduced cohomology (with $k$-coefficients) of $F_{f, \bold 0}$ is non-zero only in degree $1$, $\widetilde H^1(F_{f, \bold 0})\cong \Bbb C$, and $\widetilde T^1_{f, \bold 0}$ is the identity. Hence, $\phi_f[-1]\Bbb C_X^\bullet[2]$ is simply the extension by zero of the constant sheaf at the origin, i.e.,  $\phi_f[-1]\Bbb C_X^\bullet[2]\cong m_!\Bbb C^\bullet_\bold 0\cong m_*\Bbb C^\bullet_\bold 0$, and $\operatorname{id}-\widetilde T_{f} = 0$.

This agrees with Theorem 2.2, since 
$$IC_{{}_{V(xy)}}^\bullet\ \cong\ j^x_!\Bbb C_{V(y)}^\bullet[1]\oplus j^y_!\Bbb C_{V(x)}^\bullet[1]\ \not\cong\ \Bbb C_{V(f)}^\bullet[1].$$

There is another interesting aspect of this example. Consider the nearby cycles $\psi_f[-1]\Bbb C_X^\bullet[2]$ and their monodromy $T_f$. Then, the triple $(\operatorname{id}, T_f, \widetilde T_f)$ acts on the canonical short exact sequence of perverse sheaves
$$
0\rightarrow C_{V(f)}^\bullet[1]\rightarrow \psi_f[-1]\Bbb C_X^\bullet[2]\rightarrow \phi_f[-1]\Bbb C_X^\bullet[2]\rightarrow 0.
$$
Since $\widetilde T_f=\operatorname{id}$, it follows that $(\operatorname{id}-T_f)^2 = 0$. However, it is not difficult to prove that, in a much more general setting, the image of $\operatorname{id}-T_f$ is isomorphic to $\phi_f[-1]\Bbb C_X^\bullet[2]$. Thus, $T_f\neq \operatorname{id}$. Therefore, $T_f = \operatorname{id}+N$, where $N\neq 0$ and $N^2=0$. 

This non-trivial nilpotent portion of $T_f$ is related to a deep result of Gabber (see [{\bf Mac}], \S13). Gabber's result, applied to our current simple case, tells us that the existence of a non-trivial nilpotent portion of $T_f$, when $\widetilde T_f$ has no such non-trivial nilpotent piece, occurs precisely because $V(xy)$ is not an intersection cohomology manifold.

\vskip .2in

In Theorem 2.2, we make the above remark precise; we prove that $V(f)$ is a $k$-$IC^\bullet$ manifold if and only if there is an isomorphism of pairs of complexes and automorphisms
$$(\psi_f[-1] k_X^\bullet[n]; \ T_f)\cong (k_{{}_{V(f)}}^\bullet[n-1]\oplus\phi_f[-1]\cx; \ \operatorname{id}\oplus\widetilde T_f).$$
This result tells us that, while, in general, there may be extra structure present in the nearby cycles and their monodromy, in the case where $V(f)$ is a $k$-$IC^\bullet$ manifold, one may as well simply study the vanishing cycles and their monodromy.

\vskip .2in

Of course, the reader may suspect that if $V(f)$ has non-isolated singularities, then $V(f)$ would almost {\bf never} be a $k$-$IC^\bullet$ manifold. In fact, by combining Theorem 2.2 with the Sebastiani-Thom Theorem  ([{\bf S-T}], [{\bf O}], [{\bf Sa}],  [{\bf N1}], [{\bf N2}], [{\bf Mas4}], [{\bf Sch}]), we show in Theorem 3.2 that simple modifications of $f$ will produce a sequence of functions, $f_j$, whose zero loci are $k$-$IC^\bullet$ manifolds and whose critical loci are isomorphic to the critical locus of $f$. (Here, by ``critical locus'', we actually mean the support of the vanishing cycles.) The modifications $f_j$ are so closely related to $f$ that questions about the vanishing cycles and monodromy for arbitrary functions $f$ are reduced to the case where one may assume that $V(f)$ is an intersection cohomology manifold.

\vskip .1in

In Section 4, we give a simple ``application'' of the results of Sections 1 and 2 to line singularities. Suppose that the domain of $f$ is an open neighborhood of the origin in $\Bbb C^n$ and that $f$ has a one-dimensional critical locus through the origin. By selecting a linear form $z_0$ so that $f_{|_{V(z_0)}}$ has an isolated critical point at $\bold 0$, we may define a family of isolated singularities $f_t:=f_{|_{V(z_0-t)}}$. It is unknown, in the general case, how to effectively calculate the cohomology of the Milnor fiber, $F_{f, \bold 0}$, of $f$ at the origin. (However,  many very nice results are known; in particular, the work of Barlet in [{\bf B}] is very closely related to the results of this this paper.)

In Theorem 4.1, we recall that the reduced cohomology $\widetilde H^{n-2}(F_{f, \bold 0})$ is a monodromy-invariant subspace of both $\widetilde H^{n-2}(F_{f_0, \bold 0})$ and $\bigoplus_\nu \widetilde H^{n-2}(F_{f_{t_0}, \bold p_\nu})$, where the summation is over all irreducible components $\nu$ of $\Sigma f$ through the origin and $\bold p_\nu$ is any point in $\ob\cap\nu\cap V(z_0-t_0)$ for $0<|t_0|\ll\epsilon\ll 1$. As a corollary of Theorems 2.2 and 4.1, we conclude that if $V(f_0)$ is a $k$-intersection cohomology manifold and $V(f_{t_0})$ is not, then $\widetilde H^{n-2}(F_{f, \bold 0})$ is a {\bf proper} subspace of $\bigoplus_\nu \widetilde H^{n-2}(F_{f_{t_0}, \bold p_\nu})$.

As an example of this last corollary, consider the Whitney umbrella, presented as nodes degenerating to a cusp, given by $f = z_2^2-z_1^3-z_0z_1^2$.  As $f_{t_0}$ has a node at the origin, the shifted constant sheaf on $V(f_{t_0})$ is not intersection cohomology; intersection cohomology would ``count the origin twice''. However,  the shifted constant sheaf on $V(f_0)$ is intersection cohomology, since the cusp is a topological manifold. Therefore, $\widetilde H^{1}(F_{f, \bold 0})$ is a proper subspace of $\bigoplus_\nu \widetilde H^{n-2}(F_{f_{t_0}, \bold p_\nu})\cong  \widetilde H^{1}(F_{f_{t_0}, (t_0, \bold 0)})\cong k$, and our corollary ``explains'' the well-known fact that  $\widetilde H^{1}(F_{f, \bold 0})=0$.

\vskip .2in

We should remark here that our results are closely related to the work of many people from 15, 20, or even 40 years ago.  The results that we prove in Theorems 1.1, 2.2, and 3.2 are very closely related to work of Milnor, Sebastiani and Thom, Goresky and MacPherson,  Borho and MacPherson, Randell, Barlet, Kashiwara, Malgrange, Sabbah, and, no doubt, others. By combining old results of these mathematicians, we could arrive at many of our conclusions fairly quickly. Moreover, recent work of Torrelli also yields some of our results. However, the proofs that we give here are extremely short and do, of course, use background material that was developed during the 1980's and early 1990's. Moreover, we need the statements of our theorems in the form and generality that we give.

\vskip .2in

We would like to thank Mark Goresky, Alberto Verjovsky, Tristan Torrelli, and Daniel Barlet for making a number of helpful remarks.

\vskip .3in

\noindent\S1. {\bf Intersection Cohomology Manifolds}  

\vskip .1in

Fundamental references for many of the tools used throughout this paper are [{\bf G-M1}],  [{\bf G-M2}],  [{\bf G-M3}], and  [{\bf K-S}]. 

Throughout this paper,  we let $X$ be an $n$-dimensional complex analytic space, and let $f:X\rightarrow\Bbb C$ be a complex analytic function. By intersection cohomology, we shall always mean middle perversity intersection cohomology.

Fix a field $k$. Frequently one wishes to take $k=\Bbb Q$ or $k=\Bbb C$; however, to detect torsion, one may wish to take the finite fields $k=\Bbb Z/p\Bbb Z$. All cohomology in this paper is calculated with $k$-coefficients. We use field coefficients for two reasons: to make certain duality statements true, and so that the category of perverse complexes of sheaves of $k$-modules (vector spaces) on $X$, $Perv(X)$, is a locally Artinian category in which the simple objects are extensions by zero of intersection cohomology complexes on irreducible subvarieties with coefficients in an irreducible local system.

\vskip .3in

Below, it is convenient to adopt the convention that the empty set is a $(-1)-dimensional$ $k$-cohomology sphere. We let $\Cal D$ denote the Verdier dualizing functor on the derived category $D^b_c(X)$.

The following theorem, over the rational numbers, is due to Borho and MacPherson [{\bf B-M}]; the proof over arbitrary fields remains the same, and all of the equivalences below are either in the statement or the proof of Proposition 1.4 of [{\bf B-M}]. We give our own proof of the theorem for the sake of self-containment.

\vskip .3in

\noindent{\bf Theorem 1.1} ([{\bf B-M}]). {\it  The following are equivalent:

\vskip .1in

\noindent {\rm 1)} $X$ is a $k$-intersection cohomology manifold (i.e.,  $\ic\cong \cx$);

\vskip .1in

\noindent {\rm 2)} $\Cal D\left(\cx\right)\cong \cx$;

\vskip .1in

\noindent  {\rm 3)} for all $x\in X$, there is an isomorphism of stalk cohomology  $H^*\left(\Cal D\left(\cx\right)\right)_x\cong H^*\left(\cx\right)_x$;

\vskip .1in

\noindent  {\rm 4)} $X$ is a $k$-cohomology (or homology) manifold of real dimension $2n$, i.e., for all $x\in X$, for all $i$, $H^i(X, X-\{x\};\ k)=0$ unless $i=2n$, and $H^{2n}(X, X-\{x\};\ k)\cong k$;

\vskip .1in

\noindent  {\rm 5)} for all $x\in X$, the real link of $X$ at $x$ is a $(2n-1)$-dimensional $k$-cohomology sphere;

\vskip .1in

\noindent  {\rm 6)} there exists a Whitney stratification, $\Cal S$, of $X$ such that, for all $S\in\Cal S$, if $d$ is the complex dimension of $S$, then the real link of $S$ is a $(2(n-d)-1)$-dimensional $k$-cohomology sphere;

\vskip .1in

\noindent  {\rm 7)} for all Whitney stratifications, $\Cal S$, of $X$, for all $S\in\Cal S$, if $d$ is the complex dimension of $S$, then the real link of $S$ is a $(2(n-d)-1)$-dimensional $k$-cohomology sphere;
}

\vskip .2in

\noindent{\it Proof}. Intersection cohomology is characterized by two conditions: the support and cosupport conditions (see [{\bf G-M2}]). The support condition is always satisfied by $\cx$ on an $n$-dimensional space. In addition, the support and cosupport conditions are dual to each other when using field coefficients. Moreover, the support condition is a condition on stalk cohomology. Hence, 1), 2), and 3) are equivalent.

Now,  $H^i\left(\Cal D\left(\cx\right)\right)_x\ \cong\   H^{-i+n}(X, X-\{x\})\ \cong\ H_{-i+n}(X, X-\{x\})$, where the last isomorphism holds since our coefficients are in a field. Thus, 3) and 4) are equivalent, whether we use homology manifolds or cohomology manifolds.

Locally embed $X$, at $x$, into affine space and let $B_\epsilon(x)$ denote a small closed ball, of radius $\epsilon$, centered at $x$, in the ambient space. Then, by excision, 
$$H^{-i+n}(X, X-\{x\})\cong H^{-i+n}(B_\epsilon(x)\cap X, B_\epsilon(x)\cap X-\{x\})\cong$$
$$\widetilde H^{-i+n-1}(B_\epsilon(x)\cap X-\{x\})\cong \widetilde H^{-i+n-1}(K_{{}_{X, x}}),$$
where $K_{{}_{X, x}}$ denotes the real link of $X$ at $x$. Thus, 4) and 5) are equivalent. 

That 6) and 7) are equivalent follows from the equivalence of 1) and 5), and the fact that a complex is the intersection cohomology complex if and only if the shifted restriction of the complex to normal slices to Whitney strata yields intersection cohomology.
\qed

\vskip .3in

\noindent{\it Remark 1.2}. Mark Goresky has pointed out that we do not need complex analytic spaces or middle perversity intersection cohomology for some of the implications in Theorem 1.1; if $X$ is a reasonable space (say, homeomorphic to a simplicial complex) and the shifted constant sheaf on $X$ is isomorphic to the intersection cohomology complex, with {\bf any} perversity, then $X$ is a $k$-homology manifold.

\vskip .1in

In [{\bf T}], in the case where $X$ is a local complete intersection, Tristan Torrelli gives another characterization for $X$ being a $\Bbb C$--$IC^\bullet$ manifold ; this characterization is in terms of the Bernstein functional equation. See also Remark 2.4 for the hypersurface case.

\vskip .3in

\noindent{\bf Corollary 1.3}. {\it Suppose that $X$ is an $n$-dimensional $k$-$IC^\bullet$ manifold. 

Then, every open subset of $X$ is a $k$-$IC^\bullet$ manifold,  $X$ is purely $n$-dimensional, and the irreducible components of $X$ are disjoint from each other.

Thus, for each irreducible component $X_0$ of $X$, $k_{{}_{X_0}}^\bullet[n]$ is a simple object in $\operatorname{Perv}(X_0)$.}

\vskip .2in

\noindent{\it Proof}. Conditions 3), 5), 6), and 7) of Theorem 1.1 are clearly local; hence, every open subset of a $k$-$IC^\bullet$ manifold is a $k$-$IC^\bullet$ manifold.

If $X$ had a component, $X_0$, of dimension less than $n$, then at a generic (smooth) point of $X_0$, condition 5)  of Theorem 1.1 would fail.

Let $X_1$ be an irreducible component of $X$ such that $X_1\cap\overline{X-X_1}\neq\emptyset$. Then there is a Whitney stratification of $X$ such that $X_1\cap\overline{X-X_1}$ is a union of strata. Let $S$ be a top-dimensional stratum of $X_1\cap\overline{X-X_1}$. Then, the real link of $S$ has more than one connected component; this would contradict 7) from Theorem 1.1.
\qed

\vskip .3in

\noindent{\bf Corollary 1.4}. {\it Let $X$ and $Y$ be arbitrary analytic spaces. Then, $X\times Y$ is a $k$-$IC^\bullet$ manifold if and only if $X$ and $Y$ are both $k$-$IC^\bullet$ manifolds.}

\vskip .2in

\noindent{\it Proof}.  We believe that it is well-known that the product of homology manifolds is a homology manifold; however, for lack of a convenient reference, we supply a proof. The technical results that we contained in Corollary 2.0.4 of [{\bf Sch}].

Let $n$ and $m$ denote the (global) dimensions of $X$ and $Y$, respectively. Let $\pi_X$ and $\pi_Y$ denote the projections from $X\times Y$ onto $X$ and $Y$, respectively. Then,
$$\Cal D (k_{{}_{X\times Y}}^\bullet[n+m])\cong\pi_X^*\Cal D(k_X^\bullet[n])\overset L\to\otimes \pi_Y^*\Cal D(k_Y^\bullet[m]),$$
and so, for all $(p,q)\in X\times Y$, for all $i$,
$$
H^i\big(\Cal D (k_{{}_{X\times Y}}^\bullet[n+m])\big)_{(p, q)}\cong\bigoplus_{s+t=i}H^s\big(\Cal D(k_X^\bullet[n])\big)_p\otimes H^t\big(\Cal D(k_Y^\bullet[m])\big)_q.
$$

Using 3) from Theorem 1.1, it follows immediately that if $X$ and $Y$ are $k$-$IC^\bullet$ manifolds, then so is $X\times Y$. 

Now suppose that $X\times Y$ is a $k$-$IC^\bullet$ manifold. By selecting $p$ to be a smooth point on an $n$-dimensional component of $X$, and allowing $q$ to vary over $Y$, we find (again using 3) from Theorem 1.1) that $Y$ must be a $k$-$IC^\bullet$ manifold; the analogous argument implies that $X$ must be a  $k$-$IC^\bullet$ manifold.
\qed

\vskip .3in

\noindent{\it Remark 1.5}. Since being a $k$-intersection cohomology manifold is a local property, Corollary 1.4 immediately implies that an analytic space which fibers locally trivially over a $k$-$IC^\bullet$ manifold, with a fiber which is also a $k$-$IC^\bullet$ manifold, must itself be a $k$-$IC^\bullet$ manifold.

\vskip .3in

\noindent\S2. {\bf Hypersurfaces and Vanishing Cycles}  

\vskip .1in

Recall that $f:X\rightarrow\Bbb C$ is a complex analytic function. Let $j$ denote the closed inclusion of $V(f)$ into $X$, and let $i$ denote the open inclusion of $X-V(f)$ into $X$. We must recall a number of basic results in the derived category $D^b_c(X)$.

The composed functors $Ri_*i^*$ and $i_!i^!$ take perverse sheaves to perverse sheaves. The shifted nearby cycle functor $\psi_f[-1]$ is a functor from $D^b_c(X)$ to $D^b_c(V(f))$ which takes perverse sheaves to perverse sheaves. The nearby cycle monodromy $T_f$ is a natural automorphism of the functor $\psi_f[-1]$. The shifted vanishing cycle functor $\phi_f[-1]$ is a functor from $D^b_c(X)$ to $D^b_c(V(f))$ which takes perverse sheaves to perverse sheaves.  The vanishing cycle monodromy $\widetilde T_f$ is a natural automorphism of the functor $\phi_f[-1]$.

When the complex of sheaves, $\Adot$, under consideration is clear, we normally continue to write $T_f$ and $\widetilde T_f$ for the morphisms $T_f(\Adot):\psi_f[-1]\Adot\rightarrow\psi_f[-1]\Adot$ and $\widetilde T_f(\Adot):\phi_f[-1]\Adot\rightarrow\phi_f[-1]\Adot$.

\vskip .1in

Let $\Adot\in D^b_c(X)$.
There are four canonical distinguished triangles that we shall need:
$$
\Adot\rightarrow Ri_*i^*\Adot\rightarrow j_!j^![1]\Adot@>[1]>>\ ;\tag{i}
$$
$$
j_*j^*[-1]\Adot\rightarrow i_!i^!\Adot\rightarrow \Adot@>[1]>>\ ;\tag{ii}
$$
$$
j^*[-1]\Adot@>\operatorname{comp}>>\psi_f[-1]\Adot@>\operatorname{can}>>\phi_f[-1]\Adot@>[1]>>\ ;\tag{iii}
$$
and
$$
\phi_f[-1]\Adot@>\operatorname{var}>>\psi_f[-1]\Adot@>\operatorname{pmoc}>>j^![1]\Adot@>[1]>>\ .\tag{iv}
$$

There are equalities $\operatorname{var}\circ\operatorname{can}=\operatorname{id}-T_f$ and $\operatorname{can}\circ\operatorname{var}=\operatorname{id}-\widetilde T_f$. We let $\omega_f$  denote the morphism from $j^*[-1]\Adot$ to $j^![1]\Adot$ given by $\operatorname{pmoc}\circ\operatorname{comp}$. The octahedral lemma implies that the mapping cones of $\omega_f$ and $\operatorname{id}-\widetilde T_f$ are isomorphic; in particular, $\omega_f$ is an isomorphism if and only if $\operatorname{id}-\widetilde T_f$ is an isomorphism. We refer to $\omega_f$ as the {\it Wang morphism}.

The triples $(\operatorname{id}, T_f, \widetilde T_f)$ and $( \widetilde T_f, T_f, \operatorname{id})$ yield automorphisms of the distinguished triangles in (iii) and (iv), respectively. In particular, if $f$ is identically zero on $X$, then $\psi_f[-1]\Adot=0$ and so $\phi_f[-1]\Adot\cong j^*\Adot\cong j^!\Adot$, and $\widetilde T_f$ is the identity.

Suppose now that $\Adot$ is a simple object in $\operatorname{Perv}(X)$. If $\operatorname{supp}\Adot\not\subseteq V(f)$, then the distinguished triangles (i) and (ii) imply that $j_!j^![1]\Adot$ and $j_*j^*[-1]\Adot$ are in $\operatorname{Perv}(X)$ and, thus, $j^![1]\Adot$ and $j^*[-1]\Adot$ are in $\operatorname{Perv}(V(f))$. 

\vskip .3in

\noindent{\bf Lemma 2.1}. {\it Suppose that $X$ is an $n$-dimensional $k$-$IC^\bullet$ manifold. 

If $f$ does not vanish identically on any irreducible component of $X$, then $j^*[-1]\cx$ and $j^![1]\cx$ are perverse sheaves on $V(f)$, which are Verdier dual to each other, and the image (in $Perv(V(f))$) of $\operatorname{id}-T_f:\psi_f[-1]\cx\rightarrow\psi_f[-1]\cx$ is isomorphic to $\phi_f[-1]\cx$.

If $f$ vanishes identically on an irreducible component $X_0$ of $X$, then, when restricted to $X_0$, $\operatorname{id}-\widetilde T_f: \phi_f[-1]\cx\rightarrow\phi_f[-1]\cx$ is the zero morphism and, hence, does not induce an isomorphism in any stalk at a point of $X_0$. 
}

\vskip .2in

\noindent{\it Proof}. Suppose that $f$ does not vanish identically on any irreducible component of $X$. As the restriction of $\cx$ to each irreducible component of $X$ is a simple object,  the perversity of $j^*[-1]\cx$ and $j^![1]\cx$ follows from the case discussed above of a simple perverse sheaf. As $X$ is a $k$-$IC^\bullet$ manifold, $\Cal D(\cx)\cong\cx$, and so $\Cal D(j^*[-1]\cx)\cong j^![1]\Cal D(\cx)\cong j^![1]\cx$. Since $j^*[-1]\cx$ and $j^![1]\cx$ are perverse sheaves, the distinguished triangles (iii) and (iv) above become short exact sequences in $Perv(V(f))$ when we replace $\Adot$ with $\cx$. Thus, the morphism $\operatorname{can}$ is a surjection and the morphism $\operatorname{var}$ is an injection. It follows that the image of $\operatorname{id}-T_f = \operatorname{var}\circ\operatorname{can}$ is isomorphic to $\phi_f[-1]\cx$. This proves the first statement.

\vskip .1in

Suppose that $f$ vanishes identically on an irreducible component $X_0$ of $X$. From the paragraphs preceding the lemma, we know that, when restricted to $X_0$,   $\phi_f[-1]\cx\cong \cx$ and that $\operatorname{id}-\widetilde T_f$ is the zero morphism. This proves the second statement.\qed

\vskip .5in

Now we come to our main theorem, which relates the vanishing cycle monodromy to $k$-intersection cohomology manifolds.

\vskip .3in

\noindent{\bf Theorem 2.2}. {\it Let $X$ be an $n$-dimensional $k$-$IC^\bullet$ manifold. The following are equivalent:

\vskip .1in

\noindent {\rm a)} $\operatorname{id}-\widetilde T_f: \phi_f[-1]\cx\rightarrow\phi_f[-1]\cx$ is an isomorphism;

\vskip .2in

\noindent{\rm b)} $V(f)$ is a $k$-$IC^\bullet$ manifold, and $\operatorname{id}-\widetilde T_f: \phi_f[-1]\cx\rightarrow\phi_f[-1]\cx$ is an isomorphism when restricted to a generic subset of $V(f)$;

\vskip .2in

\noindent {\rm c)}  the pairs of complexes and morphisms $(\psi_f[-1]\cx; \ T_f)$ and $(k_{{}_{V(f)}}^\bullet[n-1]\oplus\phi_f[-1]\cx; \ \operatorname{id}\oplus\widetilde T_f)$ are isomorphic in $D^b_c(V(f))$, i.e., there is an isomorphism $\theta:\psi_f[-1]\cx\rightarrow k_{{}_{V(f)}}^\bullet[n-1]\oplus\phi_f[-1]\cx$ such that $\theta\circ T_f = \big(\operatorname{id}\oplus\widetilde T_f\big)\circ\theta$;

\vskip .1in

\noindent and, in these equivalent cases, $V(f)$ contains no component of $X$, $V(f)$ is purely $(n-1)$-dimensional, and the isomorphism in c) is an isomorphism in $Perv(V(f))$.
}

\vskip .2in

\noindent{\it Proof}. By Lemma 2.1, a) and b) both imply that  $V(f)$ contains no component of $X$. If $\psi_f[-1]\cx\cong k_{{}_{V(f)}}^\bullet[n-1]\oplus\phi_f[-1]\cx$, then, since $\psi_f[-1]\cx$ and $\phi_f[-1]\cx$ are perverse, it follows that $k_{{}_{V(f)}}^\bullet[n-1]$ is perverse; as $\cx$ is perverse when restricted to each irreducible component, this implies that $V(f)$ cannot contain a component of $X$. Therefore, by Lemma 2.1, a), b), and c) all imply that $V(f)$ contains no component of $X$, $V(f)$ is purely $(n-1)$-dimensional, and  $j^*[-1]\cx$ and $j^![1]\cx$ are perverse sheaves on $V(f)$, which are Verdier dual to each other. In particular, if the isomorphism in c) holds, then it is an isomorphism in $Perv(V(f))$.

\vskip .1in

\noindent Proof of a)$\implies$b): Suppose that  $\operatorname{id}-\widetilde T_f$ is an isomorphism. 

As $\operatorname{id}-\widetilde T_f$ is an isomorphism, the Wang morphism $\omega_f: j^*[-1]\cx\rightarrow j^![1]\cx$ is an isomorphism. Thus, 
$$k^\bullet_{{}_{V(f)}}[n-1]\cong j^*[-1]\cx\cong j^![1]\cx \cong \Cal D( j^*[-1]\cx)\cong \Cal D(k^\bullet_{{}_{V(f)}}[n-1]).$$ By 2) of Theorem 1.1, $V(f)$ is an $(n-1)$-dimensional $k$-$IC^\bullet$ manifold.

\vskip .2in

\noindent Proof of b)$\implies$a): Suppose that  $\operatorname{id}-\widetilde T_f: \phi_f[-1]\cx\rightarrow\phi_f[-1]\cx$ is an isomorphism when restricted to a generic subset of $V(f)$, and $V(f)$ is a $k$-$IC^\bullet$ manifold.

Since $\operatorname{id}-\widetilde T_f$ is an isomorphism on a generic subset of $V(f)$, the Wang morphism $\omega_f: j^*[-1]\cx\rightarrow j^![1]\cx$ is an isomorphism on a generic subset of $V(f)$. As $V(f)$ is an $(n-1)$-dimensional $k$-$IC^\bullet$ manifold,  $k^\bullet_{{}_{V(f)}}[n-1]\cong j^*[-1]\cx$  is simple when restricted to each irreducible component of $V(f)$. As $j^![1]\cx\cong\Cal D\big( j^*[-1]\cx\big)$,  $j^![1]\cx$ is also simple when restricted to each irreducible component of $V(f)$. Thus, $\omega_f$ must be an isomorphism on all of $V(f)$, and, hence, so is $\operatorname{id}-\widetilde T_f$.

\vskip .2in

\noindent Proof of a)$\implies$c): Assume a). As $j^*[-1]\cx$ is perverse, the canonical distinguished triangle (iii) from the beginning of the section becomes a short exact sequence in the category of perverse sheaves on $V(f)$ 
$$
0\rightarrow k_{{}_{V(f)}}^\bullet[n-1]@>{\operatorname{comp}}>>\psi_f[-1]\cx@>{\operatorname{can}}>>\phi_f[-1]\cx\rightarrow 0.\tag{$\dagger$}
$$
Recall that the variation natural transformation yields a morphism ${\operatorname{var}}: \phi_f[-1]\cu\rightarrow\psi_f[-1]\cu$, and ${\operatorname{can}}\circ{\operatorname{var}}= \operatorname{id}-\widetilde T_f$. 
Thus, ${\operatorname{var}}\circ(\operatorname{id}-\widetilde T_f)^{-1}$ yields a splitting of $(\dagger)$. Moreover, this splitting is compatible with the morphism of $(\dagger)$ given by the triple $(\operatorname{id}, T_f, \widetilde T_f)$. The isomorphism of  complexes and automorphisms in c) follows.

\vskip .2in

\noindent Proof of c)$\implies$a): Assume c). As $V(f)$ contains no component of $X$, the first statement of Lemma 2.1 implies that $\operatorname{im}(\operatorname{id}- T_f)\cong \phi_f[-1]\cx$. On the other hand, the isomorphism of the complexes and automorphisms implies that $\operatorname{im}(\operatorname{id}- T_f)\cong 0\oplus\operatorname{im}(\operatorname{id}- \widetilde T_f)$. Therefore, $\operatorname{im}(\operatorname{id}- \widetilde T_f)\cong \phi_f[-1]\cx$. As $Perv(V(f))$ is locally Artinian, it follows that $\operatorname{id}- \widetilde T_f$ is an isomorphism.
\qed

\vskip .3in

\noindent{\it Remark 2.3}. The condition in Theorem 2.2.b that $\operatorname{id}-\widetilde T_f$ is an isomorphism when restricted to a generic subset of $V(f)$ may seem like a very strong, unusual condition. In fact, this condition is what one expects to happen, for on the open set $\Cal U:= V(f)-\operatorname{supp}\big(\phi_f[-1]\cx\big)$, $\phi_f[-1]\cx$ is zero and thus $(\operatorname{id}-\widetilde T_f)_{|_{\Cal U}}$ is an isomorphism.

\vskip .1in

Consider, for example, the case where $X$ is an open subset of $\Bbb C^{n}$. If $f$ is reduced, then $\phi_f[-1]\cx$ is zero on a generic subset of $V(f)$, and so $\operatorname{id}-\widetilde T_f$ is an isomorphism on a generic subset of $V(f)$.  In this case, one can prove the equivalence of a) and b) in Theorem 2.2 via the Wang sequence, precisely as Milnor does in the isolated case in Chapter 8 of [{\bf Mi}]. 

Suppose now that $f$ is not reduced, and $f=f_1^{\alpha_1}f_2^{\alpha_2}\dots f_k^{\alpha_k}$ is the irreducible factorization of $f$ at a point $\bold p\in X$. At a generic point $\bold q_i\in V(f_i)$, up to isomorphism, we have the case where $f$ is a power of a coordinate function, i.e., $f=z_0^{\alpha_i}$. By taking a normal slice, we are reduced to considering the case where $f:\Bbb C\rightarrow\Bbb C$ is given by $f(z) = z^{\alpha_i}$. The Milnor fiber consists of $\alpha_i$ points and the monodromy $T_f: k^{\alpha_i}\rightarrow k^{\alpha_i}$ is simply a cyclic permutation of the coordinates. Let $\Delta$ denote the diagonal in $k^{\alpha_i}$. The map induced by $\operatorname{id}-\widetilde T_f$ is simply the automorphism of $k^{\alpha_i}/\Delta$ induced by $\operatorname{id}-T_f$. It is trivial to show that this map is an isomorphism if and only if $\alpha_i$ is not divisible by the characteristic of $k$. Therefore, we conclude that if the characteristic of $k$ does not divide the multiplicity of any component of $f$, then $\operatorname{id}-\widetilde T_f$ is an isomorphism on a generic subset of $V(f)$. 

\vskip .3in

\noindent{\it Remark 2.4}. In the case where $X$ is smooth, $f$ is reduced, and $k=\Bbb C$, Tristan Torrelli [{\bf T}] uses the language and techniques of $\Cal D$-modules to establish the same equivalence as we obtain in a) and b) of our Theorem 2.2. Moreover, Torrelli includes one more condition which is equivalent to $V(f)$ being a $\Bbb C$--$IC^\bullet$ manifold: namely, that the reduced Bernstein polynomial of $f$ has no integral root.

\vskip .5in

The following topological result was suggested to us by Alberto Verjovsky. It says that in the classical affine situation, if the singular set of the hypersurface $V(f)$ is not too big, then we have a fairly obvious generalization of the results of Milnor in Chapter 8 of [{\bf Mi}].

\vskip .3in

\noindent{\bf Theorem 2.5}. {\it Let $X$ be an open subset of $\Bbb C^n$. Let $s$ denote the dimension of the singular set $\Sigma V(f)$, and suppose that $n-1-s\geqslant 3$. Then, the following are equivalent:

\vskip .1in

\noindent {\rm 1)} $V(f)$ is an integral homology (or cohomology) manifold;

\vskip .1in

\noindent {\rm 2)} $V(f)$ is a $k$-$IC^\bullet$ manifold for $k=\Bbb Q$ (or $\Bbb C$) and for $k=\Bbb Z/p\Bbb Z$ for all primes $p$;

\vskip .1in

\noindent {\rm 3)} for all complex analytic Whitney stratifications of $V(f)$, the real links of all strata are topological spheres;

\vskip .1in

\noindent {\rm 4)} there exists a complex analytic Whitney stratification of $V(f)$ such that the real links of all strata are topological spheres;

\vskip .1in

\noindent {\rm 5)} $V(f)$ is a topological manifold.
}

\vskip .2in

\noindent{\it Proof}. For all $x\in X$, the real link of $X$ at $x$ has the homotopy-type of a finite CW-complex and, hence, the integral homology groups are finitely-generated. Therefore, the equivalence of 1) and 2) follow from the Universal Coefficient Theorems, together with Theorem 1.1.

\vskip .1in

Certainly, 3) implies 4), and 5) implies 1).

\vskip .1in

Suppose that $S$ is a Whitney stratum of $X$. Then, along $S$, $X$ is locally a product of Euclidean space with the cone on the real link of $S$. Hence, 4) implies 5).

\vskip .1in

We need to prove that 1) implies 3).

\vskip .1in

Assume that $V(f)$ is an integral homology manifold, and let $\Cal S$ be a Whitney stratification of $V(f)$. By the equivalence of 1) and 2), combined with Theorem 1.1, we find that the real link of each stratum of $\Cal S$ is an integral homology sphere. 

The real link of a stratum is the real link of a point in the normal slice. As $n-1-s\geqslant 3$, when we take the real link of a point in a normal slice to a singular stratum, that link will be at least $5$-dimensional and will be simply-connected by Theorem 5.2 of [{\bf Mi}].

Therefore, the real links of the strata are at least $5$-dimensional and are homotopy-equivalent to spheres. If we can show that the real links to strata are topological manifolds, then we will be finished by the Generalized PoincarŽ Conjecture (Theorem).

We prove this by contradiction. Let $S\in\Cal S$ be a stratum of maximum dimension such that its real link, $K$, is not a topological manifold. The link $K$ is the transverse intersection of a small sphere in the ambient space with a normal slice to $V(f)$. This link inherits a Whitney stratification $\Cal S^\prime$ from $\Cal S$. The links of strata of $\Cal S^\prime$ are precisely the links of the corresponding strata of $\Cal S$; by our maximality assumption on $S$, the links of the strata of $\Cal S^\prime$ are topological spheres. The topology of $K$ along a stratum of $\Cal S^\prime$ is a product of Euclidean space with the cone on a sphere. Thus, $K$ is a manifold; a contradiction.
\qed

\vskip .3in

\noindent\S3. {\bf Suspending $f$}  

\vskip .1in

For $j\in\Bbb N$, $j\geqslant 2$, let $f_j:X\times\Bbb C\rightarrow \Bbb C$ be given by $f_j(\bold z, w):=f(\bold z)+w^j$. We refer to $f_j$ as the {\it $(j-1)$-fold suspension of $f$}. We use this terminology because the homotopy-type of the Milnor fiber of $f_j$ at a point of $(\bold p, 0)\in V(f)\times\{0\}$ is the one-point union of $(j-1)$ copies of the suspension of the Milnor fiber of $f$ at $\bold p$; see the references below.

Our interest in suspending $f$ stems from the Sebastiani-Thom isomorphism (see [{\bf O}], [{\bf Sa}], [{\bf S-T}], [{\bf N1}], [{\bf N2}], [{\bf Sch}], [{\bf Mas4}]), which implies that $\operatorname{supp}\big(\phi_{f_j}[-1]k^\bullet_{{}_{X\times\Bbb C}}[n+1]\big) = \operatorname{supp}\big(\phi_{f}[-1]k^\bullet_{{}_{X}}[n]\big)\times\{0\}$ and, at a point 
$(\bold p, 0)\in V(f)\times\{0\}$,
$$H^{i}\left(\phi_{f_j}[-1]{k^\bullet_{{}_{X\times\Bbb C}}[n+1]}\right)_{(\bold p, 0)}\ \cong\  H^{i}\left(\phi_{f}[-1]{k^\bullet_{{}_{X}}[n]}\right)_{\bold p}\otimes k^{j-1},$$
i.e.,
$$\widetilde H^{n+i}\big(F_{{}_{f_j, (\bold p, 0)}}\big)\ \cong\  \Big(\widetilde H^{n+i-1}\big(F_{f, \bold p}\big)\Big)^{j-1},$$
and that the  action of $\widetilde T_{f_j}$ on the stalk cohomology of the vanishing cycles of $f_j$   is the tensor product of the action of $\widetilde T_{f}$ on the stalk cohomology of the vanishing cycles of $f$ at $\bold p$ and the action on $k^{j-1}$ given by a cyclically permuting the coordinates of $k^j$ and modding out by the diagonal.

\vskip .4in

\noindent{\bf Definition 3.1}. A point $\bold p\in V(f)$ is {\it totally $f$-unipotent} if and only if,  at every point $\bold q\in V(f)$ near $\bold p$, and in all degrees $i$, if $H^i\big(\phi_{f}[-1]k^\bullet_{{}_{X}}[n]\big)_\bold p\neq 0$, then the characteristic polynomial, ${\operatorname{char}}_{(\widetilde T^i_f)_\bold q}(t)$, of the stalk of $\widetilde T_f$ at $\bold q$ in degree $i$ is a power of $(t-1)$.

\vskip .1in

Note that if $\bold p\in V(f) - \operatorname{supp}\big(\phi_{f}[-1]k^\bullet_{{}_{X}}[n]\big)$, then $\bold p$ is totally $f$-unipotent.

\vskip .4in

\noindent{\bf Theorem 3.2}. {\it Let $X$ be an $n$-dimensional $k$-$IC^\bullet$ manifold. Suppose that $f$ does not vanish identically on a component of $X$. Let $\bold p\in V(f)$. Let $$J:=\{j\in\Bbb N\ |\ j\geqslant 2, V(f_j)\text{ is a $k$-$IC^\bullet$ manifold in a neighborhood of }(\bold p, 0)\}.$$ Then, $J$ is an infinite set. Moreover, if $\bold p$ is  totally $f$-unipotent, then $J=\{j\in\Bbb N\ |\ j\geqslant 2\}$.

In addition, if $\bold p$ is not totally $f$-unipotent, then $\Bbb N-J$ is also an infinite set.
}

\vskip .2in

\noindent{\it Proof}. First note that Corollary 1.4 implies that $X\times \Bbb C$ is a $k$-$IC^\bullet$ manifold.

Now, as $\operatorname{supp}\big(\phi_{f_j}[-1]k^\bullet_{{}_{X\times\Bbb C}}[n+1]\big) = \operatorname{supp}\big(\phi_{f}[-1]k^\bullet_{{}_{X}}[n]\big)\times\{0\}$ and since we are assuming that  $f$ does not vanish identically on a component of $X$, the dimension of the support of $\phi_{f_j}[-1]k^\bullet_{{}_{X\times\Bbb C}}[n+1]$ is, at most, $n-1$. However, the dimension of each component of $V(f_j)$ is precisely $n$. Therefore, $\operatorname{id}-\widetilde T_{f_j}$ is an isomorphism when restricted to a generic subset of $V(f_j)$, and Theorem 2.2 implies that $V(f_j)$ is a $k$-$IC^\bullet$ manifold if and only if $\operatorname{id}-\widetilde T_{f_j}$ is an isomorphism on all of $V(f_j)$.

Since $1$ is an eigenvalue of a linear operator over $k$ if and only if $1$ is an eigenvalue of the induced linear operator over some extension of $k$, we may assume that $k$ is algebraically closed.

Fix a Whitney stratification of $X$ such that $V(f)$ is a union of strata. Define a $\bold p$-stratum to be a stratum of this Whitney stratification which contains $\bold p$ in its closure.
Let $D$ be the set of $v\in\Bbb C$ such that there exists $x\in \operatorname{supp}\phi_{f}[-1]k^\bullet_{{}_{X}}[n]$ such that $x$ is in a $\bold p$-stratum and $v$ is an eigenvalue of the stalk at $x$ of the vanishing cycle monodromy in some degree. Then, by constructibility, $D$ is a finite set, and by the Monodromy Theorem (see  [{\bf Cl}], [{\bf Gr}], [{\bf La}], [{\bf L\^e2}])  $D$ consists of roots of unity. Let $E$ the set of natural numbers $q$ such that a primitive $q$-th root of unity is contained in $D$. Then, $E$ is a finite set. 

\vskip .1in

Note that  $\bold p$ is totally $f$-unipotent if and only if $E=\{1\}$ or $E=\emptyset$.  

\vskip .1in

If $E=\emptyset$, then $\bold p\not\in\operatorname{supp}\big(\phi_{f}[-1]k^\bullet_{{}_{X}}[n]\big)$. Thus, for all $j$, $(\bold p, 0)\not\in \operatorname{supp}\big(\phi_{f_j}[-1]k^\bullet_{{}_{X\times\Bbb C}}[n+1]\big)$, and so, in a neighborhood of $(\bold p, 0)$, $\operatorname{id}-\widetilde T_{f_j}$ induces isomorphisms on stalk cohomology (since the stalk cohomology is zero). Therefore, Theorem 2.2 implies that $J=\{j\in\Bbb N\ |\ j\geqslant 2\}$.

\vskip .1in

Suppose now that $E\neq\emptyset$.

\vskip .1in

By the Sebastiani-Thom Theorem, if $j$ is relatively prime to each element of $E$, then the vanishing cycle monodromy of $f_j$ does not have $1$ as an eigenvalue in any stalk, and hence, by Theorem 2.2, $j\in J$.

On the other hand, if $\bold p$ is  not totally $f$-unipotent and $j$ is a multiple of an element of $E-\{1\}$, then $1$ will be an eigenvalue of the vanishing cycle monodromy $\widetilde T_{f_j}$ at points near $(\bold p, 0)$; thus, $V(f_j)$ is not a $k$-$IC^\bullet$ manifold.
\qed

\vskip .3in

The trick above of adding a power of $w$ in order to eliminate $1$ as an eigenvalue of the monodromy appears in [{\bf L\^e1}].

\vskip .3in

\noindent{\it Example 3.3}. Suspending functions is interesting even in the most basic case of the ordinary quadratic singularity $f:\Bbb C^{n}\rightarrow\Bbb C$ given by $f:= z_1^2+\dots+z_n^2$, which is an iterated suspension.  The vanishing cycle monodromy is given by the single map multiplication by $(-1)^{n}$ on $\widetilde H^{n-1}(F_{f, \bold 0})\cong k$. 

When $n=2$,  one can simply see the two irreducible components which prevent the shifted constant sheaf from being intersection cohomology. More generally, if $n$ is even, the previous paragraph implies that $V(f)$ is not a $k$-$IC^\bullet$ manifold, regardless of what the field $k$ is. 

\vskip .1in

If $n$ is odd, then $V(f)$ is a $k$-$IC^\bullet$ manifold, provided that $k$ does not have characteristic $2$.

\vskip .1in

If $k$ has characteristic $2$, then the ordinary quadratic singularity {\bf never} yields a $k$-$IC^\bullet$ manifold, except in the trivial case where $n=1$.

\vskip .4in

The following is an easy corollary of Theorems 2.5 and 3.2.

\vskip .3in

\noindent{\bf Corollary 3.4}. {\it Let $X$ be an open subset of $\Bbb C^n$, and suppose that $f= f(\bold z)$ is reduced at each point of $V(f)$. If $V(f)$ is an integral homology manifold, then the subset of $X\times\Bbb C^2$ given by $V(f(\bold z)+w_1^2+w_2^2)$ is a topological manifold whose (analytic) singular set is isomorphic to the singular set of $V(f)$.
}

\vskip .2in

\noindent{\it Proof}. Let $s$ denote the dimension of the singular set $\Sigma V(f)$. As $f$ is reduced, $n-1-s\geqslant 1$. Thus, $(n+2)-1-s\geqslant 3$, and we are in a position to apply Theorem 2.5 to $V(f(\bold z)+w_1^2+w_2^2)$.  

Suppose $V(f)$ is an integral homology manifold. As we saw in the proof of Theorem 2.5, this is equivalent to $V(f)$ being a $k$-$IC^\bullet$ manifold for $k=\Bbb Q$ and $k=\Bbb Z/p\Bbb Z$, for all primes $p$. By the Sebastiani-Thom Theorem, at each point in the critical locus of $f+w_1^2+w_2^2$, the eigenvalues of the vanishing cycle monodromy of $f+w_1^2+w_2^2$ are $(-1)^2=1$ times the eigenvalues of the vanishing cycle monodromy of $f$. Thus, Theorem 2.2 and Remark 2.3 imply that $V(f(\bold z)+w_1^2+w_2^2)$ is an integral homology manifold. The desired result now follows from Theorem 2.5.
\qed

\vskip .3in

The reader should compare the above ``double suspension'' result with that of Cannon in [{\bf Ca}].

\vskip .3in

\noindent\S4. {\bf One-Dimensional Critical Loci}  

\vskip .1in

For the sake of concreteness and ease of notation, we assume throughout this section that $X=\Cal U$ is a connected open neighborhood of the origin in $\Bbb C^n$, that  $f:\Cal U\rightarrow \Bbb C$ has a one-dimensional critical locus through the origin, and that $n\geqslant 3$. Then, $f$ is reduced at the origin, and the only possibly non-zero reduced cohomology  vector spaces of $F_{f, \bold 0}$ occur in degrees $n-2$ and $n-1$.  

Select a linear form $z_0$ so that $f_{|_{V(z_0)}}$ has an isolated critical point at $\bold 0$, and define a family of isolated singularities $f_t:=f_{|_{V(z_0-t)}}$ (throughout this section, we assume that our neighborhood $\Cal U$ is re-chosen as small as necessary). Fix a small $t_0\neq 0$. In this section, we will show that if $V(f_{t_0})$ is a $k$-$IC^\bullet$ manifold and $V(f_0)$ is not, then there are restrictions on $\widetilde H^{n-2}(F_{f, \bold 0})$ (and, hence, on $\widetilde H^{n-1}(F_{f, \bold 0})$). Note that $f_0$ and $f_{t_0}$ are both reduced near the origin.

\vskip .2in

First we need to recall some fairly well-known background material. 

\vskip .2in

Let $\Pdot$ denote the perverse sheaf $\left(\phi_f[-1]\cu\right)_{{|}_{\Sigma f}}$. Assume we are in a small enough neighborhood of the origin so that $\Pdot$  is constructible with respect to $\{\Sigma f-\{\bold 0\}, \{\bold 0\}\}$. Let $\nu$ be an irreducible component of $\Sigma f$. Let $\bold p_\nu$ denote a point on $\nu-\{\bold 0\}$. Then,  $\Pdot_{|_{\nu-\bold 0}}[-1]$ is a local system with stalk isomorphic to $\widetilde H^{n-2}(F_{f_{t_0}, \bold p_\nu})\cong k^{{\overset\circ\to\mu}_\nu}$, where ${\overset\circ\to\mu}_\nu$ is the ``generic Milnor number along $\nu$'' near $\bold 0$. As $\nu-\bold 0$ is homotopy-equivalent to a circle,  $\Pdot_{|_{\nu-\bold 0}}[-1]$ is characterized by a monodromy automorphism $h_\nu$ on $\widetilde H^{n-2}(F_{f_{t_0}, \bold p_\nu})$. As $\widetilde T_f$ induces an automorphism on the complex $\Pdot$, the Milnor monodromy induced on each $\widetilde H^{n-2}(F_{f_{t_0}, \bold p_\nu})$ must commute with $h_\nu$.

 Let $\hat m:\{\bold 0\}\hookrightarrow\Sigma f$ denote the closed inclusion,  and let $\hat l:\Sigma f-\{\bold 0\}\hookrightarrow \Sigma f$ denote the open inclusion. There is a canonical distinguished triangle (see (i) at the beginning of Section 2)
$$\hat m_!\hat m^!\Pdot\rightarrow\Pdot\rightarrow R\hat l_*\hat l^*\Pdot\rightarrow \hat m_!\hat m^!\Pdot[1].\tag{$\dagger$}$$ The cosupport condition on the perverse sheaf $\Pdot$ implies that $H^{-1}(\Pdot)_\bold 0 = 0$, and so the long exact sequence on stalk cohomology obtained from $(\dagger)$ yields that there is an injection $\beta$ given by the following composition
$$\beta:\ \widetilde H^{n-2}(F_{f, \bold 0})\cong H^{-1}(\Pdot)_\bold 0\hookrightarrow H^{-1}(R\hat l_*\hat l^*\Pdot)_\bold 0\cong\oplus_\nu\operatorname{ker}(\operatorname{id}-h_\nu)\hookrightarrow \bigoplus_\nu \widetilde H^{n-2}(F_{f_{t_0}, \bold p_\nu}),$$
where $\beta$ commutes with the respective Milnor monodromies. This injection was first proved by Siersma in [{\bf Si}].

\vskip .4in

\noindent{\bf Theorem 4.1}. {\it  The vector space $\widetilde H^{n-2}(F_{f, \bold 0})$ injects into both $\widetilde H^{n-2}(F_{f_0, \bold 0})$ and  $\bigoplus_\nu \widetilde H^{n-2}(F_{f_{t_0}, \bold p_\nu})$ by inclusions which commute with the actions of the respective Milnor monodromies.

If $\dm \widetilde H^{n-2}(F_{f, \bold 0})=\sum_\nu\overset\circ\to\mu_\nu$, then each $h_v$ is the identity, and the eigenspaces of the Milnor monodromy on $\bigoplus_\nu \widetilde H^{n-2}(F_{f_{t_0}, \bold p_\nu})$ inject into  the corresponding eigenspaces of the Milnor monodromy on  $\widetilde H^{n-2}(F_{f_0, \bold 0})$.}

\vskip .2in

\noindent{\it Proof}. The discussion above, especially the definition of the injection $\beta$, tells us that  $\widetilde H^{n-2}(F_{f, \bold 0})$ injects into  $\bigoplus_\nu \widetilde H^{n-2}(F_{f_{t_0}, \bold p_\nu})$ by an inclusion which commutes with the actions of the respective Milnor monodromies. 

The Main Theorem of [{\bf Mas2}] implies that $\widetilde H^{n-2}(F_{f, \bold 0})$ injects into $\widetilde H^{n-2}(F_{f_0, \bold 0})$  by  an inclusion which commutes with the actions of the respective Milnor monodromies.

If $\dm \widetilde H^{n-2}(F_{f, \bold 0})=\sum_\nu\overset\circ\to\mu_\nu$, then $\beta$ must be an isomorphism. Hence, the second statement follows immediately from the first statement.\qed

\vskip .4in

\noindent{\bf Corollary 4.2}. {\it If  $V(f_{t_0})$ is not a $k$-$IC^\bullet$ manifold, and $V(f_{0})$ is, then  $\dm \widetilde H^{n-2}(F_{f, \bold 0})<\sum_\nu\overset\circ\to\mu_\nu$.
}

\vskip .2in

\noindent{\it Proof}.  If  $V(f_{t_0})$ is not a $k$-$IC^\bullet$ manifold, then, by Theorem 2.2, the Milnor monodromy on at least one of the $\widetilde H^{n-2}(F_{f_{t_0}, \bold p_\nu})$ must have $1$ as an eigenvalue. However, if $V(f_{0})$ is a $k$-$IC^\bullet$ manifold, then Theorem 2.2 implies that $1$ is not an eigenvalue of the Milnor monodromy on $\widetilde H^{n-2}(F_{f_0, \bold 0})$. It follows at once from Theorem 4.1 that $\dm \widetilde H^{n-2}(F_{f, \bold 0})<\sum_\nu\overset\circ\to\mu_\nu$.\qed

\vskip .3in

\noindent\S5. {\bf Questions and Remarks}  

\vskip .1in

Corollary 4.2 is an interesting little result. It leads one to ask the following more general question:

\vskip .1in

\noindent{\bf Question}. If $f$ has a one-dimensional critical locus at the origin, and $\dm \widetilde H^{n-2}(F_{f, \bold 0})=\sum_\nu\overset\circ\to\mu_\nu$, does it follow that the critical locus is itself smooth and that $f$ defines a family of isolated singularities with constant Milnor number?

\vskip .1in

In fact, the example of $f:=(y^2-x^3)^2+w^2$ shows us that, in general, the answer to the above question is ``no''. However, if all of the components of the critical locus of $f$ are smooth at the origin, then we have showed in a recent paper [{\bf L-M}], with L\^e D. T., that the answer to the above question is ``yes''.

\vskip .2in

We were trying to answer the above question when we proved most of the results of this paper. We hoped that, by suspending $f$, we could reduce ourselves to the case where $V(f)$ is an intersection cohomology manifold, and that we could then make effective use of the nice categorical property of intersection cohomology: namely,  that the simple objects in the locally Artinian category of perverse sheaves are precisely the intersection cohomology complexes on irreducible subvarieties with coefficients in irreducible local systems. 

\vskip .2in

One might also hope to gain more information by considering the nearby cycles $\psi_f[-1]\cu$ instead of the vanishing cycles $\phi_f[-1]\cu$. Theorem 2.2 tells us that the nearby cycles may contain interesting structure that is not present in the vanishing cycles precisely when $V(f)$ is {\bf not} an intersection cohomology manifold. Thus, it may be useful to suspend $f$ in such a way that we do {\bf not} have that $V(f)$ is is an intersection cohomology manifold; this is one of the reasons that we look at the set $\Bbb N-J$ is Theorem 3.2.

Consider, for instance,  a basic example which shows that the nearby cycles are ``more interesting'' than the vanishing cycles when $V(f)$ is not an intersection cohomology manifold; let $f:\Bbb C^2\rightarrow\Bbb C$ be given by $f(x, y) = x^2+y^2$. As we saw in Example 3.3, the monodromy map  $\widetilde T_f$ on the vanishing cycles is simply the identity in this case. Does this imply that the monodromy map $T_f$ is also the identity map? No. Recall that, in $\operatorname{Perv}(V(f))$, the image of $\operatorname{id}-T_f$ is isomorphic to the vanishing cycles. Thus, if the vanishing cycles are non-zero, then $T_f$ is not the identity. However, it is true that if $\widetilde T_f =\operatorname{id}$, then $(\operatorname{id}-T_f)^2 = 0$.

\enddocument